\documentclass[11pt]{amsart}
\usepackage{amsfonts}
\usepackage{amssymb}
\usepackage{latexsym} 
\usepackage{graphicx} 
\usepackage{pifont} 

\newtheorem{theorem}{Theorem}[section]
\newtheorem{definition}[theorem]{Definition}
\newtheorem{lemma}[theorem]{Lemma}
\newtheorem{corollary}[theorem]{Corollary}

\newtheorem{example}[theorem]{Example}
\newtheorem{question}[theorem]{Question}
\newtheorem{remark}[theorem]{Remark}
\newtheorem{acronym}[theorem]{Acronym}

\newcommand\NNN{{\mathbb N}}

\newcommand\ZZZ{{\mathbb Z}}

\newcommand\iv{^{-1}} 

\newcommand\cat{^{\mathord{\frown}}} 
\newcommand\nil{\,()\,} 

\renewcommand\iff{\longleftrightarrow}
\newcommand\imp{\longrightarrow} 

\newcommand{\pref}[1]{(\ref{#1})} 

\newcommand\sbl[1]{\langle#1\rangle} 
\newcommand\ld{\backslash} 
\newcommand\rd{/} 
\newcommand{\Mlt}{\mathrm{Mlt}} 

\begin{document}

\title[A Generalization of Moufang
and Steiner Loops]{A Generalization of Moufang and Steiner Loops}

\author[M.K. Kinyon]{Michael K. Kinyon}
\address{Department of Mathematics\\
Western Michigan University\\
Kalamazoo, MI\ 49008-5248 USA}
\email{mkinyon@wmich.edu}
\urladdr{http://unix.cc.wmich.edu/\symbol{126}mkinyon}
\author[K. Kunen]{Kenneth Kunen$^*$}
\thanks{$^*$Author supported by NSF Grant DMS-9704520.}
\address{Department of Mathematics\\
University of Wisconsin\\
Madison, WI 57306 USA}
\email{kunen@math.wisc.edu}
\urladdr{http://www.math.wisc.edu/\symbol{126}kunen}
\author{J.D. Phillips}
\address{Department of Mathematics \& Computer Science\\
Wabash College\\
Crawfordsville, IN\ 47933 USA}
\email{phillipj@wabash.edu}
\urladdr{http://www.wabash.edu/depart/math/faculty.htm{\#}Phillips}
\thanks{This paper is in final form, and no version of it will be submitted for
publication elsewhere.}
\date{\today}

\begin{abstract}
We study a variety of loops, RIF, which arise naturally from
considering inner mapping groups, and a somewhat larger
variety, ARIF. All Steiner and Moufang loops are RIF, and all
flexible C-loops are ARIF. All ARIF loops are diassociative.
\end{abstract}

\maketitle

\section{Introduction}
\label{sec-intro}

A \textit{loop} is an algebraic system $(L;\; \cdot, \ld, \rd, 1)$
satisfying the equations
\[
x \cdot (x \ld y) = x \ld (x \cdot y) =
(y \rd x) \cdot x = (y \cdot x) \rd x = 
y\cdot 1 = 1\cdot y = y \ \ .
\]
See the books \cite{BR, CHPS, PF} for further information. Since
loops in general form too broad a class for detailed study, the
literature has focused on various sub-varieties of loops.

Many of these varieties are defined by some weakening of the
\textit{associative} law, $x \cdot yz = xy \cdot z$. Some obvious
weakenings are the \textit{flexible} laws and the \textit{left}
and \textit{right} \textit{alternative} laws:
\[
\mathrm{FLEX:} \ x \cdot yx = xy \cdot x \quad
\mathrm{RALT:} \ x \cdot yy = xy \cdot y \quad
\mathrm{LALT:} \ y \cdot yx = yy \cdot x .
\]
There is also the \textit{inverse property, IP}. This asserts that
there is a permutation $J$ of order two such that (writing $x\iv$
for $xJ$) left and right division are given by $x \ld y = x\iv y$
and $y \rd x = y x\iv$. To see that the IP defines a variety of
loops, note that it can be expressed by the equations
$x\ld y = (x\rd 1)\, y$ and $y\rd x = y\, (1\rd x)$. Most of the loops
considered in this paper have the IP. The IP implies the
antiautomorphic inverse property (AAIP), $(xy)\iv = y\iv x\iv$,
so that $J$ provides an isomorphism from the loop $(L; \cdot)$
onto its opposite loop $(L; \circ)$ (where
$x \circ y = y \cdot x$). Thus, in IP loops, the right and left
versions of properties (e.g., RALT and LALT) are equivalent.

In a loop $L$, the left and right translations by $x\in L$ are
defined by $yL(x)=xy$ and $yR(x)=yx$, respectively. The
\textit{multiplication group} of $L$ is the permutation group on
$L$, $\Mlt(L) = \langle R(x),L(x) : x\in L\rangle$, generated by
all left and right translations. The \textit{inner mapping group}
is the subgroup $\Mlt_1(L)$ fixing $1$. If $L$ is a group, then
$\Mlt_1(L)$ is the group of inner automorphisms of $L$. In an IP
loop, the AAIP implies that we can conjugate by $J$ to get:
\[
L(x)^J = R(x^{-1}) \quad R(x)^J = L(x^{-1})
\]
where $\theta^J = J\iv \theta J = J\theta J$ for a permutation
$\theta$. If $\theta$ is an inner mapping, then so is $\theta^J$.
This leads us to one of the classes of IP loops we study in this
paper:

\begin{definition}
A \emph{RIF loop} is an IP loop $L$ with the property that
$\theta^J = \theta$ for all $\theta\in \Mlt_1(L)$. Equivalently,
inner mappings preserve inverses, i.e.,
$(x\iv)\theta = (x\theta)\iv$ for all $\theta\in \Mlt_1(L)$ and
all $x\in L$.
\end{definition}

RIF loops include the \textit{Steiner loops\/}, which may be
defined to be IP loops of exponent two (that is, $x\iv = x$, so
$J$ is the identity permutation). Steiner loops arise naturally
in combinatorics, since they correspond uniquely to Steiner
triple systems; specifically, the Steiner loop $L$ corresponds
to the triple system
$\{\{x, y ,xy\} : x \ne y \ \&\ x,y \ne 1\}$ on $L \setminus \{1\}$.

RIF loops also include what is probably the most well-studied class
of nonassociative loops, namely those satisfying the
\textit{Moufang laws} \cite{M1, M2}:

\begin{definition}\label{def:moufang}
A \emph{Moufang loop} is a loop satisfying the following equations:
\[
\begin{array}{rlrl}
M1:& (x(yz))x = (xy)(zx)\qquad\qquad &
M2:& (xz)(yx) = x((zy)x) \\
N1:& ((xy)z)y = x(y(zy))\qquad\qquad &
N2:& ((yz)y)x = y(z(yx))
\end{array}
\]
\end{definition}
In fact, by work of Bol and Bruck, each of these equations implies
the other three (see Bruck \cite{BR}, Lemma 3.1, p.~115). That
every Moufang loop is RIF follows from Lemma 3.2, p.~117, of 
\cite{BR}. It is easily seen that the only loops which are both
Steiner and Moufang are the boolean groups. Thus a direct product
of a nonassociative Steiner loop with a nonassociative Moufang
loop is a RIF loop which is neither Steiner nor Moufang.

A loop is said to be \textit{diassociative} if the subloop
$\sbl{x,y}$ generated by any two elements is a group. Diassociative
loops are always IP loops, and are flexible and alternative.
Steiner loops are obviously diassociative; in fact each $\sbl{x,y}$
is a boolean group (of order $1$, $2$, or $4$). Less obviously, by
Moufang's Theorem, every Moufang loop is diassociative.

Bruck and Paige \cite{BP} defined an \textit{A-loop} to be a loop
in which every inner mapping is an automorphism. An A-loop need
not be an IP-loop, but they show, by modifying the proof of Moufang's
Theorem, that every IP A-loop is diassociative. In fact, it turned
out later \cite{KKP} that the IP A-loops form a proper sub-variety
of the Moufang loops. Weakening of the notion of IP A-loop so that
inner mappings preserve inverses, but not necessarily products, we
obtain RIF loops.

The notion ``RIF'' can be expressed by a finite set of equations
(Lemma \ref{lem:RIF-equiv}). These lead naturally (Lemma
\ref{lem:RIF-implies-ARIF}) to a slightly weaker notion, ARIF.

\begin{definition}
\label{def:ARIF}
An ARIF loop is a flexible loop satisfying the following equations:
\[
W1: R(x) R(yxy) = R(xyx) R(y) \qquad
W2: L(x) L(yxy) = L(xyx) L(y)
\]
\end{definition}
In fact, a flexible loop satisfying either W1 or W2 has the IP and hence
satisfies both equations (Lemma \ref{lem:ARIF-IP}). Every ARIF loop
of odd order is Moufang (Corollary \ref{cor-odd}) (whereas non-group
Steiner loops are RIF and not Moufang). Besides RIF loops, ARIF
includes another variety of IP loops, namely the flexible C-loops
(Corollary \ref{coro:flex-C-ARIF}).
C-loops were introduced by Fenyves \cite{FEN};
see Section \ref{sec-basics}. There
exist flexible C-loops which are not RIF loops (Example
\ref{ex-flex-C-non-RIF}), and there exist ARIF loops which are neither RIF
loops nor C-loops (see Section \ref{sec-ex}).

\begin{acronym}
$A = \mathit{Almost}$,
$R = \mathit{Respects}$,
$I = \mathit{Inverses}$,
$F = \mathit{Flexible}$.
\end{acronym}

Section \ref{sec-dias} is devoted to the proof of our main result,
a generalization of Moufang's Theorem to ARIF loops:

\begin{theorem}
\label{thm-dias}
Every ARIF loop is diassociative.
\end{theorem}

Our inductive proof of this theorem is patterned on Moufang's
proof, but is quite a bit more complicated than hers, or than
the corresponding proof in Bruck and Paige \cite{BP} for IP
A-loops. We do not know a simpler proof, but Example
\ref{ex-steiner} shows that the basic lemma on associators
developed by Moufang can fail in a ARIF loop (in fact, in a
Steiner loop).

Note that if we write out the definition of diassociativity
in the obvious way, we get an infinite list of equations.
The following problem, asked first by Evans and Neumann
\cite{ER}, is still open:

\begin{question}
\label{question-dias}
Does the variety of diassociative loops have a finite basis?
\end{question}
If the answer is ``yes'', which seems unlikely, then inductive
proofs of diassociativity could always be replaced by the
verification of a finite number of instances of diassociativity,
which could result in a simplification.

Our investigations were aided by the automated reasoning tools
OTTER, developed by McCune \cite{MC}, and SEM developed by
J. Zhang and H. Zhang \cite{ZZ}. SEM finds finite models of
systems of axioms, and was used to produce the three examples
in Section \ref{sec-ex}. OTTER derives statements from axioms,
and was used to derive enough instances of diassociativity
from ARIF for the pattern to become clear.

\section{Basics}
\label{sec-basics}

Following Bruck \cite{BR} (see IV.1), the inner mapping group
of any loop is generated by the inner mappings of the form
$L(x,y)$, $R(x,y)$, and $T(x)$:

\begin{definition}\label{def:inner}
\qquad $T(x) = R(x)L(x)\iv$ \\
$L(x,y) = L(x)L(y)L(yx)\iv$ \qquad
$R(x,y) = R(x)R(y)R(xy)\iv$ 
\end{definition}

Using this, we can express the notion of RIF by equations.

\begin{lemma}
\label{lem:RIF-equiv}
The following are equivalent for an IP loop $L$:
\begin{itemize}
\item[1.] $L$ is a RIF loop.
\item[2.] $L$ is flexible and $R(x,y) = L(x\iv, y\iv)$ for
all $x,y\in L$.
\item[3.] $R(xy)L(xy) = L(y)L(x)R(x)R(y)$ for all $x,y\in L$.
\item[4.] $L(xy)R(xy) = R(x)R(y)L(y)L(x)$ for all $x,y\in L$.
\end{itemize}
\end{lemma}
\begin{proof}
The flexible law can be expressed as $R(x)L(x) = L(x)R(x)$ for all
$x$. In an IP loop, this is equivalent to
$L(x\iv)R(x) = R(x)L(x\iv)$, that is, $T(x)^J = T(x)$. Also, an easy
calculation gives $R(x,y)^J = L(x\iv, y\iv)$ in an IP loop. Thus
(1) and (2) are equivalent. Using the IP and Definition
\ref{def:inner}, $R(x,y) = L(x\iv, y\iv)$ is equivalent to
$L(xy)R(xy) = L(y)L(x)R(x)R(y)$. Since the flexible law is just
$R(z)L(z) = L(z)R(z)$, (2) implies (3). Conversely, if (3) holds,
then taking $y=1$ gives the flexible law, so that (3) implies (2).
Finally, (3) and (4) are equivalent by the IP.
\end{proof}

Combining 3 and 4 from Lemma \ref{lem:RIF-equiv} we obtain the very useful
identity $L(y)L(x)R(x)R(y) = R(x)R(y)L(y)L(x)$, which we will frequently 
appeal to in our arguments.

For the next result, we introduce the notation $P(x) = L(x)R(x)$.

\begin{corollary}
In a RIF loop, $P(xyx) = P(x)P(y)P(x)$.
\end{corollary}
\begin{proof}
Applying Lemma \ref{lem:RIF-equiv} twice,
$P(x\cdot yx) = 
R(x)R(yx)L(yx)L(x) = 
R(x)\; L(x)L(y)R(y)R(x)\; L(x) = 
P(x)P(y)P(x) $.
\end{proof}

The fact that Moufang loops satisfy $P(xyx) = P(x)P(y)P(x)$ is
Theorem 5.1, p.~120, of Bruck \cite{BR}. The same theorem points
out that $L(xyx) = L(x)L(y)L(x)$ and $R(xyx) = R(x)R(y)R(x)$ also
hold in Moufang loops. But in flexible loops, these are simply
restatements of the Moufang equations $N1,N2$ in Definition
\ref{def:moufang}, so they do not hold in all RIF loops, since
they fail in any non-group Steiner loop.

Next we show that RIF loops satisfy equations W1, W2 of Definition
\ref{def:ARIF}.

\begin{lemma}
\label{lem:RIF-implies-ARIF}
Every RIF loop is an ARIF loop.
\end{lemma}
\begin{proof}
Equations W1, W2 are equivalent in IP loops. To prove W1, start
with $R(v)R(y)L(y)L(v) = L(y)L(v)R(v)R(y)$,
which is
\[
v (y (z v) y) = (v (y z) v) y\ \ ,
\]
and set $v = ux$ and $z = u\iv$, so that $zv = x$. We get
\[
ux \cdot yxy = (ux \cdot yu\iv \cdot ux) y\ \ .
\]
But
$R(u\iv) R(ux) L(ux) = R(x)L(x)L(u)$
(see Lemma \ref{lem:RIF-equiv}), so
\[
ux \cdot yxy = (u \cdot xyx) y \ \ ,
\]
which is W1.
\end{proof}

Next we show that every flexible C-loop is an ARIF loop.
\textit{C-loops}, introduced by Fenyves \cite{FEN}, are loops
satisfying the equation $((x y) y) z = x (y (y z))$. These have
the IP (see \cite{FEN}, Theorem 4) and are alternative (see
\cite{FEN}, Theorem 3). They are not necessarily flexible (see
Example \ref{ex-C-non-flex}). Every Steiner loop is trivially a
C-loop; in fact, Table 1 of \cite{FEN}, a C-loop which is not
Moufang, is just the 10-element Steiner loop.

\begin{theorem}
Every C-loop satisfies
\[
R(xy)^2 = R(x) R(y(xy)) = R((xy)x) R(y).
\]
\end{theorem}
\begin{proof}
Since the loop is alternative, the C-loop property can be written
as $R(a)^2 R(b) = R(a^2 b)$. This gives us:
\[
R(xy)^2 R(y\iv) = R((xy)^2 y\iv) = R((xy)((xy)y\iv)) = R((xy)x) \ \ .
\]
so $ R(xy)^2 = R((xy)x) R(y)$. Now, if $x = v\iv$ and $y = v(uv)$,
we have $xy = uv$ and hence $R(uv)^2 = R(u)R(v(uv))$.
\end{proof}

\begin{corollary}
\label{coro:flex-C-ARIF}
Every flexible C-loop is a ARIF loop.
\end{corollary}

We now examine basic properties of ARIF loops.

\begin{lemma}
\label{lem:ARIF-IP}
A loop satisfying
\[
W1': \qquad R(x)R((yx)y) = R(x(yx))R(y) \qquad
\]
is an alternative IP loop.
\end{lemma}
\begin{proof}
Let $x\iv$ denote $1 \rd x$, so that $x\iv x=1$. Applying
W1$^{\prime}$ to $x\iv$ gives $(yx)y = (x\iv (x(yx)))y$, and
cancelling gives $yx = x\iv (x(yx))$. Replacing $y$ with
$(x\ld y)\rd x$ yields $x\ld y = x\iv y$. In particular
$1 = x(x\iv 1) = xx\iv$, and so $(x\iv)\iv = x$. Next apply
W1$^{\prime}$ to $(x(yx))\iv$ to get $((x(yx))\iv x)((yx)y) = y$,
and thus $(yx)y = ((x(yx))\iv x)\iv y$. Cancelling yields
$yx = ((x(yx))\iv x)\iv$. Replacing $y$ with $y\rd x$ gives
$y = ((xy)\iv x)\iv$ and so $y\iv = (xy)\iv x$, which implies
$(xy)y\iv = x$. Replacing $x$ with $x\rd y$ gives $xy\iv = x\rd y$.
Thus the loop satisfies the IP. Setting $y = 1$ in W1$^{\prime}$
yields the right alternative law $R(x)R(x) = R(xx)$, and the right
and left alternative laws are equivalent in IP loops.
\end{proof}

\begin{corollary}
\label{coro:ARIF-IP}
Every ARIF loop is an alternative IP loop.
\end{corollary}

\begin{lemma}
\label{lemma-premoufang}
Every ARIF loop satisfies $R(x) R(y^2 x\iv) R(x) = R(xy^2)$
and $L(x) L(x\iv y^2) L(x) = L(y^2 x)$.
\end{lemma}
\begin{proof}
The second equation is equivalent to the first in IP loops.
So start with $R(aba)R(b) = R(a)R(bab)$. \\
Set $b = xy$ and $a = x\iv$ (so $ab = y$) 
to get
\[
R(yx\iv)R(xy) = R(x\iv)R(xy^2).
\]
Set $b = x$ and $a = yx\iv$ (so $ab = y)$
to get
\[
R(y^2 x\iv)R(x) = R(yx\iv)R(xy).
\]
Putting these together, we have
\[
R(y^2 x\iv)R(x) = R(x\iv)R(xy^2).
\]
\end{proof}

\begin{corollary}
\label{cor-remoufang}
Every ARIF loop in which each element is a square is a Moufang loop.
\end{corollary}
\begin{proof}
Now we have $R(x)R(y x\iv) R(x) = R(xy)$.
If we let $z = y x\iv$ and $y = zx$, we get
$R(x)R(z) R(x) = R(xzx)$, which (in flexible loops) is
the Moufang equation N1 of Definition \ref{def:moufang}.
\end{proof}

In general, products of an element of a loop with
itself do not associate, e.g., $x\cdot xx \neq xx\cdot x$.
We shall see that this problem does not arise in ARIF loops.
Until then, we let $x^n$ denote the right-associated product.

\begin{definition}
Define $x^n = (1)(L(x))^n$ for any $n \in \ZZZ$.
\end{definition}

Thus, $x^3 = x\cdot xx$, and (in an IP loop)
$x^{-3} = (1)L(x\iv)^{3} = x\iv \cdot x\iv x\iv$,
whereas $(x^3)\iv = x\iv x\iv \cdot x\iv$.

\begin{definition}
A loop $L$ is \emph{power associative} iff for all $x\in L$,
the subloop $\sbl{x}$ generated by $x$ is a group, and
\emph{power alternative} iff
$L(x^i) = (L(x))^i$ and
$R(x^i) = (R(x))^i$ 
for all $x \in L$ and all $i \in \ZZZ$
\end{definition}

It is easily seen that diasociativity implies power alternativity
and power alternativity implies power associativity.

Now for $n > 0$, let us say that an IP loop $L$ is
\textit{$n$--PA} iff $L(x^m) = (L(x))^m$ whenever $1 \le m \le n$
and $x \in L$. So, the $1$--PA is trivial and the $2$--PA (that is,
$xx\cdot y = x\cdot xy$) is equivalent to the alternative law.
Hence, a $2$--PA loop satisfies $x^3 = xx\cdot x = x\cdot xx$ and
$x^{-3} = (x^3)\iv$. An IP loop which is $n$--PA for all $n > 0$
is power alternative. Note that the $n$--PA implies that
$x^i \cdot x^j = (1 L(x)^j) L(x)^i = x^{i + j}$ whenever
$1 \leq i,j \leq n$. For $0 < j < k$, set $m = jk$, and
suppose that $L$ is $k$--PA. Then
$x^m = 1(L(x)^j)^k = (1)L(x^j)^k = (x^j)^k$ by the $j$--PA, 
so $L(x^m) = L(x^j)^k = L(x)^m$ by the $k$--PA and $j$--PA. Thus
$L$ is $m$--PA, so the smallest $n$ such that the
$n$--PA fails must be prime.

\begin{theorem}
Every ARIF loop is power alternative.
\end{theorem}
\begin{proof}
By Corollary \ref{coro:ARIF-IP} and the preceding remarks, it
is sufficient to prove that the ARIF loop $L$ is $n$--PA for
all odd $n \geq 3$. So, for $n = 2k+1$, $k\geq 1$, assume that
$L$ is $2k$--PA. Setting $y = x^k$ in Lemma
\ref{lemma-premoufang}, we get
$L(x) L(x\iv (x^k)^2) L(x) = L((x^k)^2) x)$. Now
$(x^k)^2 = x^{2k}$ and
$x\iv x^{2k} = x\iv (x x^{2k-1}) = x^{2k-1}$, and so, by the
$(2k-1)$--PA, we have $L(x)^{2k+1} = L(x^{2k} x)$. Applying this
to $1$ gives $x^{2k+1} = x^{2k} x$, so  $L(x)^n = L(x^n)$,
which is the $n$--PA.
\end{proof}

\begin{corollary}
\label{cor-odd}
Every finite ARIF loop of odd order is Moufang.
\end{corollary}
\begin{proof}
In a power alternative loop $L$, the subloop $\sbl{x}$ generated by
a given $x\in L$ induces a coset decomposition of $L$, and so if
$L$ is finite, the order of $x$ must divide the order of $L$. Thus
in a power alternative loop of odd order, each element is a square.
Now apply Corollary \ref{cor-remoufang}.
\end{proof}

\section{Diassociativity}
\label{sec-dias}

Moufang loops are diassociative by Moufang's Theorem.
The same holds for ARIF loops
(Theorem \ref{thm-dias}), as we show in this section.
First, a lemma which generalizes Lemma \ref{lemma-premoufang}:

\begin{lemma}
\label{lemma-rrll}
In any ARIF loop:
\begin{itemize}
\item[1.] $R(y x^{m}) R(x^{n} y\iv) = R(y x^{m+k}) R(x^{n-k} y\iv)$
\item[2.] $R(x^{m}y ) R(y\iv x^{n}) = R(x^{m+k}y ) R(y\iv x^{n-k})$
\item[3.] $L(x^{m}y ) L(y\iv x^{n}) = L(x^{m+k}y ) L(y\iv x^{n-k})$
\item[4.] $L(y x^{m}) L(x^{n} y\iv) = L(y x^{m+k}) L(x^{n-k} y\iv)$
\end{itemize}
whenever $m,n,k \in \ZZZ$ and either $k$ is even or $m+n$ is even.
\end{lemma}
\begin{proof}
We focus on (1,2), since (3,4) are equivalent by the IP. Let $L$
be a ARIF loop. On $\ZZZ^2$, define three relations by
\begin{eqnarray*}
(m,n)\sim_1 (s,t) &\iff& 
\forall x,y\in L \, [
R(y x^{m}) R(x^{-n} y\iv) = R(y x^{s}) R(x^{-t} y\iv) ] \\
(m,n)\sim_2 (s,t) &\iff& 
\forall x,y\in L \, [
R(x^{-m}y ) R(y\iv x^{n}) = R(x^{-s}y ) R(y\iv x^{t}) ] \\
(m,n)\sim (s,t) &\iff& (m,n)\sim_1 (s,t)\mbox{ and } (m,n)\sim_2 (s,t) \ \ .
\end{eqnarray*}
Now $(m,n)\sim_1 (s,t) \iff (m,s)\sim_2 (n,t)$, and so
$$
(m,n)\sim (s,t) \iff (m,s)\sim (n,t) \ \ .
\eqno{(A)}
$$
It is clear that each of $\sim_1,\sim_2, \sim$ is an equivalence relation.
By (A) and the fact that $\sim$ is symmetric:
$$
(m,n)\sim (s,t) \iff (n,m)\sim (t,s) \ \ .
\eqno{(B)}
$$
Also, replacing $x$ by $x\iv$ we have
$$
(m,n)\sim (s,t) \iff (-m,-n)\sim (-s,-t) \ \ .
\eqno{(C)}
$$
Replacing $y$ by $yx^j$ we have
$$
(m,n)\sim (s,t) \iff (m+j,n+j)\sim (s+j,t+j) \ \ .
\eqno{(D)}
$$
So far, everything we have said holds in any IP power
alternative loop.
Our goal is now to prove $(m,n) \sim (m+k,n+k)$
whenever $m,n,k \in \ZZZ$ and either $k$ is even or $m+n$ is even.

In the equations
\[
R(yxy)R(x) = R(y)R(xyx)
\ \ ;\ \ R(xyx)R(y) = R(x)R(yxy) \ \ ,
\]
set $x = a^\alpha b\iv$ and $y = ba^\delta$. Then,
by power alternativity,
\[
xy = a^{\alpha + \delta} a^{-\delta}  b\iv \, \cdot\,  ba^\delta =
a^{\alpha + \delta}  (ba^\delta)\iv \, \cdot\,   ba^\delta =
a^{\alpha + \delta}\ \ ,
\]
so that $xyx = a^{2\alpha + \delta} b\iv$ and
$yxy = b a^{\alpha + 2\delta}$.
We get:
\begin{eqnarray*}
R( b a^{\alpha + 2\delta} )R( a^\alpha b\iv ) 
&=& R( ba^\delta )R( a^{2\alpha + \delta} b\iv ) \\
R( a^{2\alpha + \delta} b\iv)R( ba^\delta )
&=& R( a^\alpha b\iv )R( b a^{\alpha + 2\delta} ) \ \ .
\end{eqnarray*}
The first of these equations implies
$(\alpha + 2\delta, -\alpha) \sim_1 (\delta, -2\alpha-\delta)$, while
the second
implies $(-2\alpha-\delta, \delta) \sim_2 (-\alpha, \alpha + 2\delta)$,
so (B) yields
$ (\alpha + 2\delta, -\alpha) \sim (\delta, -2\alpha-\delta) $.
Set $\alpha = -m-2c$ and $\delta = m+c$
to get $(m, m+2c) \sim (m+c, m+3c)$. Iterating this:
$$
(m, m+2c) \sim (m+jc, m+(j+2)c) \eqno{(E)}
$$
for every $m,c,j \in \ZZZ$.
But by (A), we also have $(m, m+c) \sim (m+2c, m+3c)$, and iterating this
we get:
$$
(m, m+c) \sim (m+2jc, m+(2j+1)c) \eqno{(F)}
$$
for every $m,c,j \in \ZZZ$.

Now, in view of (D), the lemma is equivalent to:
$$
n \mbox{ even or } k \mbox{ even } \imp (0,n) \sim (k,n+k) \ \ .
\eqno{(*)}
$$
We prove by induction on $n$ that $(*)$ holds for all $k$.
By (C), it is sufficient to consider $n \ge 0$, and the $n = 0$ case holds by
the IP. Now, fix $n > 0$.

If $n$ is even, we need to prove $(0,n) \sim (k,n+k)$ for all $k$.
Setting $c = \frac{n}{2}, m = k$ in (E) we get
$(k, n+k) \sim (k+j\frac{n}{2}, n + k+j\frac{n}{2})$, so it is sufficient
to prove $(0,n) \sim (k,n+k)$ whenever $0 \le k < \frac{n}{2}$.
But since this is the same as
$(0,k) \sim (n,n+k)$, it follows by applying $(*)$ inductively
to $k$, since $n$ is even.

If $n$ is odd, we need to prove $(0,n) \sim (2k,n+2k)$ for all $k$.
Setting $c = n, m = 2k$ in (F) we get
$(2k, n + 2k) \sim (2k+2jn, n + 2k+2jn) $, so it is sufficient to prove
$(0,n) \sim (2k,n+2k)$ whenever $0 \le 2k < 2n$.
Now $n$ is odd, so $2k \ne n$. If
$0 \le 2k < n$, then 
$(0,n) \sim (2k,n+2k)$ (equivalently, $(0,2k) \sim (n,n+2k)$)
follows by applying $(*)$ inductively to $2k$.
If $n < 2k < 2n$, then induction gives us instead
$(0,2n - 2k) \sim (-n,n -2k)$, and hence (by (A,C))
$(0,n) \sim (2k - 2n ,2k -n)$.
But also $(0,n) \sim (-2n,-n)$ (by (F) with $c = n, m = 0, j = -1$),
so $(2k,n + 2k) \sim (2k-2n,2k-n)$ (by (D)), and hence
$(2k,n + 2k) \sim (0,n)$.
\end{proof}

We remark that one cannot remove the restriction on $m,n,k$.
For example, if 
$R(y x) R(y\iv) = R(y) R(x y\iv)$
(that is $(1,0) \sim (0, -1)$) holds, then the loop must be Moufang
(see the proof of Corollary \ref{cor-remoufang}).
Conversely, Moufang loops satisfy the lemma
for all $m,n,k$. To see this, note that we now have
$(m+ 1,m) \sim (m, m-1)$ for every $m$, and hence
$(m+ 1,m) \sim (n + 1, n)$ for every $m,n$.
So, $(m,n) \sim (m+1,n+1)$ for every $m,n$,
and hence $(m,n) \sim (m+k,n+k)$ for every $m,n,k$.

The following lemma will be useful in the proof of diassociativity:

\begin{lemma}
\label{lemma-powers}
In a ARIF loop, suppose that $p,a,q$ satisfy:
\begin{eqnarray*}
&&p \cdot a q = p a \cdot q \\
&&p a \cdot a\iv q = p a\iv \cdot a q = p q \ \ .
\end{eqnarray*}
Then $p a^m \cdot a^n q = p a^{m + k} \cdot a^{n -k} q$ for all $m,n,k$.
\end{lemma}
\begin{proof}
We first verify
\[
p \cdot a\iv q = p a\iv \cdot q \ \ :
\]
Applying Definition \ref{def:ARIF} twice,
$R(x) R(y) R(xyx) = R(xyx) R(y) R(x)$.
Let $x = q$ and $y = q\iv a\iv$, so $xyx = a\iv q$.
Let $z = pa$. Then
\[
z R(x) R(y) = (pa \cdot q)( q\iv a\iv ) =(p \cdot aq)( q\iv a\iv ) = p \ \ ,
\]
so $z R(x) R(y) R(xyx) = p \cdot a\iv q$.
Also, $z R(xyx) = pa \cdot a\iv q = pq$, so
\[
z R(xyx) R(y) = pq \cdot q\iv a\iv = 
( p a\iv \cdot a q ) \cdot q\iv a\iv = p a\iv \ \ ,
\]
so $z R(xyx) R(y) R(x) = p a\iv \cdot q$.

{ 
\newcommand\same[4] {p a^{#1} \cdot a^{#2} q = p a^{#3} \cdot a^{#4} q}
\newcommand\RR[4]{ R(q\iv a^{#1})R(a^{#2} q) = R(q\iv a^{#3}) R(a^{#4} q)}
\newcommand\LL[4]{ L(a^{#1} p\iv)L(p a^{#2}) = L(a^{#3} p\iv ) L(p a^{#4})}
Apply $\RR {-1} {s+1} 0 s $ to $\same 0 1 1 0$ to get
$\same 0 {s+1} 1 s$ whenever $s$ is even.
Then apply $\LL 0 t {-1} {t+1}$ to get
$\same t {s+1} {t+1} s$ whenever $s,t$ are even.
Now, the same argument starting from $\same 0 {-1} {-1} 0$ results in
$\same {t} {s+1} {t-1} {s+2}$
whenever $s,t$ are even. Applying these with
$(s,t), (s+2, t-2), (s+4, t-4), \cdots$, we get
$\same {t + i} {s + 1 -i} {t+j} {s + 1 -j}$ for all $i,j$
whenever $s,t$ are even. 
But this implies that $\same m n {m+k} {n-k}$ whenever $m+n$ is odd.

Now apply $\RR {-1} {s+1} 0 s $ to $\same {-1} 1 0 0$ to get
$\same {-1} {s+1} 0 s$ whenever $s$ is even.
Then $\LL 1 {t-1} 0 {t}$ yields
$\same {t-1} {s+1} {t} s$ whenever $s,t$ are even.
The same argument starting from $\same 1 {-1} 0 0$ results in
$\same {t+1} {s-1} {t} s$ whenever $s,t$ are even.
Iterating as before, we get
$\same {t+i} {s-i} {t + j} {s -j}$ for all $i,j$ whenever $s,t$ are even,
which implies $\same m n {m+k} {n-k}$ whenever $m+n$ is even.
} 
\end{proof}

A special case of this lemma is where $p,q$ are both powers of some
element $b$. Now, in a flexible power alternative loop,
\[
x^i (y x^j) = (y)R(x)^j L(x)^i = 
(y)L(x)^i R(x)^j = (x^i y) x^j\ \ ,
\]
so the notation
$x^i y x^j$ is unambiguous.

\begin{lemma}
\label{lemma-3assoc}
In a ARIF loop, 
$x^i y^m \cdot y^n x^j = x^i y^{m+n} x^j$ for all $i,j,m,n \in \ZZZ$.
\end{lemma}
\begin{proof}
We apply Lemma \ref{lemma-powers}.
$x^i \cdot y x^j = x^i y \cdot x^j$ holds by power alternativity.
But also
\[
x^i y\iv \cdot y x^j = 
(x^i y\iv )(y x^{-i} \cdot x^{i+j} ) = x^{i+j} 
\]
by the IP, and likewise $ x^i y \cdot y\iv x^j $.
\end{proof}

\begin{remark}
A commutative flexible power alternative loop which satisfies
$(x^i y^m )(y^n x^j) = x^i y^{m+n} x^j$ for all $i,j,m,n \in \ZZZ$
is diassociative.
\end{remark}
\begin{proof}
Fix $a,b$, and let $L = \{(a^i b^m ) : i,m \in \ZZZ\}$.
By commutativity,
$(a^i b^m )(a^j b^n) = a^{i+j} b^{m+n}$,
which implies both that $L$ is a subloop and that $L$ is associative.
\end{proof}

In particular, every commutative ARIF loop is diassociative.
To prove diassociativity in the non-commutative case, we set up
some machinery. For a set $A$, let $A^*$ be the set of finite
sequences (or words) from $A$; i.e., $A^*$ is the free monoid
generated by $A$. For a word $W \in A^*$, let $|W|$ be the
length of $W$, so $|\nil| = 0$ and $|(a,b,c)| = 3$. In particular, 
$|\cdot | : A^* \to \NNN$ is just the unique homomorphic 
extension of $A \to \{1\}$. If $W,V \in A^*$, then $W\cat V$ will
denote their concatenation.

Now let $L$ be a loop. For $B,C\subseteq L$, let 
$B \cdot C = \{b\cdot c : b\in B \ \&\ c\in C\}$.
For a word $W \in A^*$ from some $A\subseteq L$, we wish
to define the set of products of $W$ under all possible
associations. We do this inductively as follows.

\begin{definition}
Define $\Pi(\nil) = \{1\}$ and
$\Pi( ( x ) ) = \{x\}$, and, when $|W| \ge 2$:
\[
\Pi(W) = \bigcup\{\Pi( V_1)\cdot \Pi( V_2) :
V_1 \cat V_2 = W \ \&\ 
V_1 \ne \nil\ \&\ V_2 \ne \nil\} \ \ .
\]
$W$ \emph{associates} iff $\Pi(W)$ is a singleton.
\end{definition}

\noindent Alternatively, $\Pi(W)$ can be described as follows.
For $A\subseteq L$, let $A'$ be the free groupoid-with-identity
on $A$. Let $p:A'\to A^*$ and $q:A'\to L$ be the unique
homomorphic extensions of the identity map on $A$. Then for
$W\in A^*$, $\Pi(W) = q(p^{-1}\{W\})$.

Among the products in $\Pi(W)$ is the right associated product
$\pi_R(W)$ defined inductively as follows.

\begin{definition}
Define $\pi_R(\nil) = 1$ and
$\pi_R( ( x ) ) = x$, and, when $|W| \ge 1$:
\[
\pi_R( ( x ) \cat W) = x \cdot \pi_R(W) \ \ .
\]
\end{definition}

\noindent Alternatively, for $A\subseteq L$, let $L^*:A^* \to \Mlt(L)$
be the unique homomorphic extension of
$L: A \to \Mlt(L); x\mapsto L(x)$. Then the right associated
product of $W\in A^*$ is $\pi_R(W)=1L^*(W)$. 
$W$ associates iff $\Pi(W) = \{\pi_R(W)\}$.

\begin{lemma}
A loop $L$ is diassociative iff for all $a,b\in L$, 
every $W \in \{a,b,a\iv,b\iv\}^*$ associates.
\end{lemma}

Now, one might try to prove that all such $W$ associate by
induction on $|W|$, in which case the following definition
and lemma might be helpful:

\begin{definition}
\label{def-pik}
If $W = (a_1, \ldots a_n)$ and $1 \le k \le n-1$, then
$ \pi^k(W) = \pi_R(a_1, \ldots,a_k) \cdot \pi_R(a_{k+1}, \ldots, a_n)$
\end{definition}

\begin{lemma}
\label{lemma-W-assoc}
If $W = (a_1, \ldots a_n)$, then $W$ associates iff:
\begin{itemize}
\item[1.] $\pi^k(W) = \pi^j(W)$ whenever $1 \le j,k \le n-1$, and
\item[2.] The words
$(a_1, \ldots,a_k)$ and $(a_{k+1}, \ldots, a_n)$
associate whenever $1 \le k \le n-1$.
\end{itemize}
\end{lemma}

In our proof that ARIF loops are diassociative, we induct not on
$|W|$, but on the \textit{block length} of $W$, defined as follows:

\begin{definition}
$B(\nil) = 0$ and $B((x)) = 1$. If $W = (x,y)\cat V$,
then $B(W) = B((y)\cat V)$ if $x \in \{y,y\iv\}$, and 
$B(W) = B((y)\cat V) + 1$ otherwise.
\end{definition}

Thus, $B(a,a,a\iv, b,b\iv, b, a) = 3$ if $a \ne b$ and $a \ne b\iv$.

\begin{definition}
An IP loop $L$ is $D$ -- associative iff
for all $a,b \in L$, every $W \in \{a,b,a\iv,b\iv\}^*$
such that $B(W) \le D$ associates.
\end{definition}

\begin{lemma}
For any IP loop $L$:
\begin{itemize}
\renewcommand\labelitemi{\ding{43}}
\item
$L$ is power associative iff 
$L$ is $1$ -- associative.
\item
$L$ is power alternative iff 
$L$ is $2$ -- associative.
\item
$L$ is diassociative iff $L$ is 
$D$ -- associative for all $D$.
\end{itemize}
\end{lemma}
So, we already know that every ARIF loop is $2$ -- associative.

\begin{lemma}
\label{lemma-dias-induct}
Suppose that an IP loop $L$ is $(D-1)$ -- associative, and
$D\nobreak\ge\nobreak3$.
Then $L$ is \emph{$D$ -- associative} iff
whenever $2 \le i \le D-1$,
$x,y \in L$, and
$n, k, m_1, m_2, \ldots, m_D \in \ZZZ$, the appropriate one
of the following equations holds:
\begin{eqnarray}
\lefteqn{ (x^{m_1} y^{m_2} x^{m_3} \cdots y^{m_i}) \cdot
(y^n x^{m_{i+1}} \cdots y^{m_{D}}) =} \hspace{1cm} \nonumber \\
& & (x^{m_1} y^{m_2} x^{m_3} \cdots y^{m_i - k}) \cdot
(y^{n+k} x^{m_{i+1}} \cdots y^{m_{D}}) \label{eveneven} \\
\lefteqn{ (x^{m_1} y^{m_2} x^{m_3} \cdots x^{m_i}) \cdot
(x^n y^{m_{i+1}} \cdots y^{m_{D}}) =} \hspace{1cm} \nonumber \\
& & (x^{m_1} y^{m_2} x^{m_3} \cdots x^{m_i - k}) \cdot
(x^{n+k} y^{m_{i+1}} \cdots y^{m_{D}}) \label{evenodd} \\
\lefteqn{ (x^{m_1} y^{m_2} x^{m_3} \cdots y^{m_i}) \cdot
(y^n x^{m_{i+1}} \cdots x^{m_{D}}) =} \hspace{1cm} \nonumber \\
& & (x^{m_1} y^{m_2} x^{m_3} \cdots y^{m_i - k}) \cdot
(y^{n+k} x^{m_{i+1}} \cdots x^{m_{D}}) \label{oddeven} \\
\lefteqn{ (x^{m_1} y^{m_2} x^{m_3} \cdots x^{m_i}) \cdot
(x^n y^{m_{i+1}} \cdots x^{m_{D}}) =} \hspace{1cm} \nonumber \\
& & (x^{m_1} y^{m_2} x^{m_3} \cdots x^{m_i - k}) \cdot
(x^{n+k} y^{m_{i+1}} \cdots x^{m_{D}}) \label{oddodd} 
\end{eqnarray}
$i$ is even in (\ref{eveneven},\ref{oddeven}) and odd in
(\ref{evenodd},\ref{oddodd}), and $D$ is even in
(\ref{eveneven},\ref{evenodd}) and odd in
(\ref{oddeven},\ref{oddodd}).
\end{lemma}

Note that by $(D -1)$ -- associativity, the parenthesized
expressions in Lemma \ref{lemma-dias-induct}
are unambiguous.
Also, note that by power alternativity, it is not necessary to
consider the cases $i = 1$ and $i = D$.
By Lemma \ref{lemma-3assoc},

\begin{corollary}
Every ARIF loop is $3$ -- associative.
\end{corollary}

Now, in proving $D$ -- associativity by induction on $D$,
equations (\ref{eveneven},\ref{evenodd},\ref{oddeven},\ref{oddodd})
give us four different cases to consider. Case \pref{oddodd} is
handled easily by conjugation. First, note that
$3$ -- associativity implies that conjugation commutes with powers:
\begin{lemma}
\label{lemma-conj}
In any 3 -- associative IP loop,
$ (x\iv y x)^n = x\iv y^n x$ for all $n \in \ZZZ$.
\end{lemma}
\begin{proof}
This is clear for $n = 0$ and $n = \pm 1$, so it is sufficient
to prove it for $n \ge 1$, which we do by induction on $n$.
Assume it holds for $n$. By 3 -- associativity,
$ x^{n+1} = x y \cdot y\iv x^{n}$.
Let $x = u\iv vu = u\iv v^2\cdot v\iv u$ and $y = u\iv v$.
Then $xy = u\iv v^2$ and $x^n = u\iv v^n u = u\iv v \cdot v^{n-1} u $,
so $y\iv x^n = v^{n-1} u$.
Hence, $(u\iv v u)^{n+1} = x^{n+1} = x y \cdot y\iv x^{n} =
u\iv v^2 \cdot v^{n-1} u = u\iv v^{n+1} u $.
\end{proof}

\begin{lemma}
\label{lemma-dias-induct-oddodd}
Suppose an IP loop $L$ is $(D-1)$ -- associative, where $D \ge 4$,
and assume that $2 \le i \le D-1$ and $D,i$ are both odd.
Then equation \pref{oddodd} of Lemma \ref{lemma-dias-induct} holds.
\end{lemma}
\begin{proof}
Under the substitution
$u = x^{m_1} y x^{-m_1}$, $y = x^{-m_1} u x^{m_1}$,
equation \pref{oddodd} reduces to:
\begin{eqnarray*}
\lefteqn{ ( u^{m_2} x^{m_3} \cdots x^{m_i + m_1}) \cdot
(x^{n -m_1} u^{m_{i+1}} \cdots x^{m_{D} + m_1}) =} \hspace{1cm} \\
& & ( u^{m_2} x^{m_3} \cdots x^{m_i + m_1 - k}) \cdot
(x^{n -m_1 + k} u^{m_{i+1}} \cdots x^{m_{D} + m_1}) \ \ ,
\end{eqnarray*}
which is an instance of $(D-1)$ -- associativity.
\end{proof}

Next, observe that in ARIF loops, Lemma \ref{lemma-powers} implies
that we need only consider
(\ref{eveneven},\ref{evenodd},\ref{oddeven},\ref{oddodd}) in two special
cases:

\begin{lemma}
\label{lemma-dias-special}
Suppose that a ARIF loop $L$ is $(D-1)$ -- associative, and
$D \ge 3$. Fix $i$ with $2 \le i \le D-1$,
and \emph{fix} $m_1, m_2, \ldots,m_{i-1}, m_{i+1} \ldots, m_D \in \ZZZ$. 
Fix $x,y\in L$.
Assume that the appropriate equation from
(\ref{eveneven},\ref{evenodd},\ref{oddeven},\ref{oddodd}) in
Lemma \ref{lemma-dias-induct} holds in the three special cases
$m_i = k = -n = 1$, 
$m_i = k = -n = -1$, and
$m_i = k = 1 \; ;\; n = 0$.
Then the same equation holds for all values of $m_i,k,n$.
\end{lemma}
\begin{proof}
For example, say $D$ and $i$ are even, so we are considering 
equation \pref{eveneven}. Let
$p = x^{m_1} y^{m_2} x^{m_3} \cdots x^{m_{i-1}}$ and let
$q = x^{m_{i+1}} \cdots y^{m_{D}}$. Then the three special cases
give us $p y \cdot y\iv q = pq$, 
$p y\iv \cdot y q = pq$, and
$p y \cdot q = p \cdot yq$. But then Lemma \ref{lemma-powers}
yields \pref{eveneven} for all $m_i,k,n$.
\end{proof}

Actually, we shall combine the first two cases and handle
$m_i = k = -n$ in Lemma \ref{lemma-dias-induct-cancel}.
First, a preliminary lemma, which is a variant of Lemma \ref{lemma-powers}.

\begin{lemma}
\label{lemma-central-assoc}
In a ARIF loop, suppose that $p,a,q,s$ are elements such that:
\begin{itemize}
\item[$\alpha$.] 
$p \cdot a^m s = p a^m \cdot s$
\ ; \ $s\iv \cdot a^m q = s\iv a^m \cdot q$
\ ; \ $p \cdot a^m q = p a^m \cdot q$ \ \ .
\item[$\beta$.] $s\iv a^m s \cdot s\iv q = s\iv a^m q$ \ \ .
\item[$\gamma$.] $p s \cdot s\iv a^m q = p a^m s \cdot s\iv q = p a^m q $ \ \ .
\item[$\delta$.] $p a^{m} s \cdot s\iv a^{-m} q = pq$ \ \ .
\end{itemize}
for all $m \in \ZZZ$. Then
\begin{itemize}
\item[$\Delta$.] $p a^{m} s \cdot s\iv a^{-n} q = p a^{m-n}q$
\end{itemize}
for all $m,n\in\ZZZ$.
\end{lemma}
\begin{proof}
Let $v = s\iv a s$ and $u = s\iv q$.
By ($\beta$) and Lemma \ref{lemma-conj}, 
$v^j u = s\iv a^j q$ and hence
$u\iv v^j = q\iv a^j s$ for every $j$.
By Lemma \ref{lemma-rrll}, 
$R(v^{-n} u) = R(v^{-m} u) R(u\iv v^{m + k}) R(v^{-n-k} u)$ whenever
$k$ is even or $m + n$ is even.
Applying this to $p a^{m} s $ and using ($\delta$), we have
\[
p a^{m} s \cdot s\iv a^{-n} q = 
[ pq \cdot q\iv a^{m+k} s ] \cdot s\iv a^{-n-k} q \ \ .
\]
But ($\delta$) also implies that $pq = p a^{m+k} s \cdot s\iv a^{-m-k} q$,
so by the IP we have
\[
p a^{m} s \cdot s\iv a^{-n} q = 
p a^{m+k} s \cdot s\iv a^{-n-k} q \ \ .
\]
If $k$ equals either $-m$ or $-n$, then this yields
$p a^{m} s \cdot s\iv a^{-n} q = p a^{m-n}q$ by ($\gamma$).
So, let $k = -m$ if $m$ is even and let $k = -n$ if $n$ is even.
If $m,n$ are both odd, then $m + n$ is even and there is no
restriction on $k$, so $k$ can be either $-m$ or $-n$.
\end{proof}

\begin{lemma}
\label{lemma-dias-induct-cancel}
Suppose a ARIF loop $L$ is $(D-1)$ -- associative, where $D \ge 4$,
and assume that $1 < i < D$. Then the appropriate equation
(\ref{eveneven},\ref{evenodd},\ref{oddeven},\ref{oddodd}) 
from Lemma \ref{lemma-dias-induct} holds whenever
$m_i = k = -n$.
\end{lemma}
\begin{proof}
First, consider \pref{eveneven}. When
$m_i = k = -n$, this reduces to:
\[
\arraycolsep=0pt
\begin{array}{rcl}
(x^{m_1} y^{m_2} x^{m_3} \cdots x^{m_{i-1}} y^{m_i}) &\cdot&
(y^{-m_i} x^{m_{i+1}} \cdots y^{m_{D}}) = \\
x^{m_1} y^{m_2} x^{m_3} \cdots &x&^{m_{i-1} + m_{i+1} }
\cdots y^{m_{D}} \ \ . 
\end{array}
\]
If we let $u = y^{-m_i} x y^{m_i}$ and $x = y^{m_i} u y^{-m_i}$,
then this becomes
\[
\arraycolsep=0pt
\begin{array}{rcl}
(y^{m_i} u^{m_1} y^{m_2} u^{m_3} \cdots u^{m_{i-1}}) &\cdot&
( u^{m_{i+1}} \cdots y^{m_{D}-m_i}) = \\
y^{m_i} u^{m_1} y^{m_2} u^{m_3} \cdots &u&^{m_{i-1} + m_{i+1} }
\cdots y^{m_{D}-m_i} \ \ . 
\end{array}
\]
which is an instance of $(D-1)$ -- associativity. A similar argument works
in cases \pref{evenodd} and \pref{oddodd} but not in case
\pref{oddeven}, where $D$ is odd and $i$ is even.

To illustrate case \pref{oddeven}, consider $D = 7$ and $i$ $=$ $2$,
$4$, or $6$. If $i = 6$, we must verify
\[
(x^{m_1} y^{m_2} x^{m_3} y^{m_4} x^{m_5} y^{m_6}) \cdot
(y^{-m_6} x^{m_{7}}) = 
x^{m_1} y^{m_2} x^{m_3} y^{m_4} x^{m_5 + m_{7}} \ \ .
\]
This is no problem, since it is equivalent to
\[
x^{m_1} y^{m_2} x^{m_3} y^{m_4} x^{m_5} y^{m_6} =
(x^{m_1} y^{m_2} x^{m_3} y^{m_4} x^{m_5 + m_{7}}) \cdot
(x^{-m_7} y^{m_{6}}) \ \ ,
\]
which is an instance of $6$ -- associativity. Likewise, the case
$i = 2$ is no problem. But, when $i = 4$, we must verify
\[
(x^{m_1} y^{m_2} x^{m_3} y^{m_4}) \cdot
(y^{-m_4} x^{m_{5}} y^{m_6} x^{m_{7}}) = 
x^{m_1} y^{m_2} x^{m_3 + m_{5}} y^{m_6}x^{m_{7}} \ \ .
\]
This is equivalent to
\[
x^{m_1} y^{m_2} x^{m_3} y^{m_4} =
(x^{m_1} y^{m_2} x^{m_3 + m_{5}} y^{m_6}x^{m_{7}}) \cdot
(x^{-m_7} y^{-m_{6}} x^{-m_5} y^{m_{4}}) \ \ ,
\]
which requires $8$ -- associativity. However, this equation
requires only $6$ -- associativity in the special case that 
$m_3 = -m_5$, and this case is sufficient by
Lemma \ref{lemma-central-assoc},
applied with
$s = y^{m_4}$, $a = x$, $p = x^{m_1} y^{m_2}$, and $q = y^{m_6}x^{m_{7}}$.
The special case is condition ($\delta$) of
Lemma \ref{lemma-central-assoc}, and
conditions $(\alpha,\beta,\gamma)$ are verified using $5$ -- associativity.

The general situation is handled similarly. We must verify
\begin{eqnarray*}
\lefteqn{ (x^{m_1} y^{m_2} x^{m_3} \cdots x^{m_{i-1}} y^{m_i}) \cdot
(y^{-m_i} x^{m_{i+1}} \cdots x^{m_{D}}) =} \hspace{1cm} \\
& & x^{m_1} y^{m_2} x^{m_3} \cdots x^{m_{i-1} + m_{i+1} }
\cdots x^{m_{D}} \ \ ,
\end{eqnarray*}
where $D$ is odd and $i$ is even.
By mirror symmetry, we may assume that $i \ge (D + 1) / 2$.
Fix $D,i,x,y$.
Let $H(r)$ be the assertion that this equation holds in
the special case that $m_{i+\ell} = -m_{i - \ell}$
whenever $1 \le \ell \le r$.
So, we want to show $H(0)$.
Now, $H(r)$ holds for $r$ large enough by $(D - 1)$ -- associativity,
and $H(r+1) \imp H(r)$ holds by Lemma \ref{lemma-central-assoc},
so we are done.

To be more specific, $H(r)$ asserts that
\[
\arraycolsep=0pt
\begin{array}{rcl}
(x^{m_1} y^{m_2} \cdots 
z^{m_{i-r-2}} w^{m_{i-r-1}} z^{m_{i-r}} 
&\cdots& x^{m_{i-1}} y^{m_i}) \cdot \\
\qquad (y^{-m_i} x^{-m_{i-1}} &\cdots&
z^{-m_{i-r}} w^{m_{i+r+1}} z^{m_{i+r+2}} 
\cdots x^{m_{D}}) = \\
x^{m_1} y^{m_2} \cdots
z^{m_{i-r-2}} &w&^{m_{i-r -1} + m_{i+r+1} } z^{m_{i+r+2}} 
\cdots x^{m_{D}} \ \ ,
\end{array}
\]
where $(z,w)$ is $(x,y)$ if $r$ is odd and $(y,x)$ if $r$ is even.
This is of form $d b = c$, which is equivalent to $d = c b^{-1}$.
Now, $c$ has $D -2r -2$ blocks and $b^{-1}$ has $D-i+1$ blocks,
and $c$ ends with $x$ while $b^{-1}$ begins with $x$, so that
the expression $c b^{-1}$ has $2D -2r -i -2$ blocks, 
so $H(r)$ follows from $(D - 1)$ -- associativity whenever
$2D -2r -i -2 \le D -1$, or $r \ge (D -i-1)/2$.

Now, assume that $r \le (D -i-1)/2 -1$ and assume that $H(r+1)$ holds.
$H(r+1)$ is the special case of $H(r)$ with $ m_{i+r+1} = -m_{i-r -1} $.
We conclude $H(r)$ by applying Lemma \ref{lemma-central-assoc}, with
$a = w$, $s = z^{m_{i-r}} \cdots x^{m_{i-1}} y^{m_i} $,
$p = x^{m_1} y^{m_2} \cdots z^{m_{i-r-2}}$, and
$q = z^{m_{i+r+2}} \cdots x^{m_{D}}$.
Condition ($\delta$) is $H(r+1)$, and the conclusion, ($\Delta$),
is $H(r)$.
We must verify that conditions ($\alpha,\beta,\gamma$) require
only $(D - 1)$ -- associativity. ($\alpha,\gamma$) are easy.
For ($\beta$), the expression $s\iv w^m s s\iv q$ has
no more than 
\[
(r+1) + 1 + (r+1) + (r+1) + (D -i -r -1) -1 = D -i + 2r +2
\]
blocks. Since $2r+2 \le D -i -1$ and $2i \ge D+1$, we have
\[
D -i + 2r +2 \le 2D - 2i -1 \le D - 2 \ \ .
\]
\end{proof}

By this lemma and Lemma \ref{lemma-dias-special}, the requirement
for $D$ -- associativity simplifies to Lemma \ref{lemma-dias-assoc}:

\begin{definition}
$ W(x,y; m_1, m_2, m_3, \ldots, m_D)$ denotes the word of length $D$,
$(x^{m_1}, y^{m_2}, x^{m_3}, \ldots, z^{m_D})$,
where $z$ is $x$ if $D$ is odd and $y$ if $D$ is even
\end{definition}

\begin{lemma}
\label{lemma-dias-assoc}
Suppose a ARIF loop $L$ is $(D-1)$ -- associative, where $D \ge 4$.
Then $L$ is $D$ -- associative iff
$ W(x,y; m_1, m_2, m_3, \ldots, m_D)$ associates for every
$x,y\in L$ and every 
$m_1, m_2, m_3, \ldots, m_D \in \ZZZ$.
\end{lemma}

To aid in proving this associativity:

\begin{lemma}
\label{lemma-dias-two}
Suppose a ARIF loop $L$ is $(D-1)$ -- associative, where $D \ge 4$,
and $W = W(x,y; m_1, m_2, m_3, \ldots, m_D)$.
Then $\pi^k(W) = \pi^{k+2}(W)$
(see Definition \ref{def-pik})
whenever $1 \le k \le D-3$.
\end{lemma}
\begin{proof}
Say $k$ is even; the argument for odd $k$ is the same.
Then we must prove 
\begin{eqnarray*}
(x^{m_1} y^{m_2} x^{m_3} \cdots y^{m_{k}}
x^{m_{k+1}} y^{m_{k+2}} ) \cdot (x^{m_{k+3}} \cdots z^{m_{D}}) = \\
(x^{m_1} y^{m_2} x^{m_3} \cdots y^{m_{k}}) \cdot 
( x^{m_{k+1}} y^{m_{k+2}} x^{m_{k+3}} \cdots z^{m_{D}}) 
\end{eqnarray*}
We apply Lemma \ref{lemma-powers}, with
$p = x^{m_1} y^{m_2} x^{m_3} \cdots y^{m_{k}}$,
$q = x^{m_{k+3}} \cdots z^{m_{D}}$, and
$a = y^{-m_{k+2}} x^{-m_{k+1}}$.
Now $ p \cdot a q = p a \cdot q $ follows by $(D-2)$ -- associativity,
and $p a \cdot a\iv q = p a\iv \cdot a q = p q $ follows by
$(D-1)$ -- associativity plus Lemma \ref{lemma-dias-induct-cancel}.
So, $ p \cdot a\iv q = p a\iv \cdot q $ follows by
Lemma \ref{lemma-powers}.
\end{proof}

\begin{lemma}
\label{lemma-dias-odd}
Suppose a ARIF loop $L$ is $(D-1)$ -- associative, where $D \ge 5$
and $D$ is odd. Then $L$ is $D$ -- associative.
\end{lemma}
\begin{proof}
If $W = W(x,y; m_1, \ldots, m_D)$, then
Lemma \ref{lemma-dias-induct-oddodd} implies 
$\pi^2(W) = \pi^3(W)$. Thus, applying Lemma \ref{lemma-dias-two},
the $\pi^k(W)$ (for $1 \le k \le D-1$) are all the same.
It follows by Lemma \ref{lemma-W-assoc} that 
$W$ associates, so $L$ is $D$ -- associative by
Lemma \ref{lemma-dias-assoc}.
\end{proof}

\begin{lemma}
\label{lemma-dias-even}
Suppose that a ARIF loop $L$ is $(D-1)$ -- associative, where $D \ge 4$
and $D$ is even. Then $L$ is $D$ -- associative.
\end{lemma}
\begin{proof}
Again, we must show that each $ W(x,y; m_1, \ldots, m_D)$ associates.
Let $H(r)$ be the assertion that
$W(x,y; m_1, \ldots, m_D)$ associates for all $x,y \in L$
and $m_1, \ldots m_D \in \ZZZ$ with $m_i = 1$ whenever $r < i \le D$.
So, our lemma is equivalent to $H(D)$.
Let $H^+(r)$ be the assertion that
$W(x,y; m_1, \ldots, m_D)$ associates for all $x,y \in L$
and $m_1, \ldots m_D \in \ZZZ$ with $m_i = 1$ whenever $r < i < D$.
So, our lemma is also equivalent to $H^+(D-1)$.
We shall in fact prove:
\begin{enumerate}
\item\label{basis}
$H(1)$.
\item\label{ind}
$H(k-1) \imp H(k)$ whenever $2 \le k \le D-1$.
\item\label{switch}
$H(D-1) \imp H^+(1)$.
\item\label{indplus}
$H^+(k-1) \imp H^+(k)$ whenever $2 \le k \le D-1$.
\end{enumerate}
Applying these items in order yields $H^+(D-1)$ and hence the lemma.

First, note that,
as in the proof of Lemma \ref{lemma-dias-odd},
$W = W(x,y; m_1, \ldots, m_D)$ associates
if 
$\pi^k(W) = \pi^{k+1}(W)$ for \textit{some} $k$ with $1 \le k \le D-2$.

To prove $H(1)$: 
Let $W = W(x,y; m, 1, 1, \ldots, 1)$. We prove that $W$ associates
by showing that $\pi^1(W) = \pi^2(W)$; that is,
$ x^m \cdot y x y \cdots y = x^m y \cdot x y \cdots y$.
Letting $u = xy$ so $y = x\iv u$, this reduces to
$x^m ( x\iv u \cdot u^{D-1}) = x^{m-1}u \cdot u^{D-1}$,
which is true by $2$ -- associativity.

To prove $H(k-1) \imp H(k)$ when $2 \le k \le D$ and $k$ is odd:
$W$ is now
$(x^{m_1}, y^{m_2}, x^{m_3},
\ldots, y^{m_{k-1}}, x^{m_k}, y^1,x^1 \ldots y^{1})$,
and we shall prove that $\pi^{k-1}(W) = \pi^{k}(W)$.
Let $p = x^{m_1} y^{m_2} x^{m_3} \cdots y^{m_{k-1}}$ and
$q = y^1 x^1 \cdots y^{1} = y (xy)^{(D-k-1)/2}$.
We need to show that $p\cdot x^{m_k} q = p x^{m_k} \cdot q$.
When $m_k = 1$, this is true by $H(k-1)$. But also
$p x \cdot x\iv q = p x\iv \cdot x q = p q$ by
Lemma \ref{lemma-dias-induct-cancel}.
$H(k)$ now follows by Lemma \ref{lemma-powers}.

The proofs for $H(k-1) \imp H(k)$ for $k$ odd and 
for $H^+(k-1) \imp H^+(k)$ are the same.

Finally, we assume $H(D-1)$ and prove $H^+(1)$.
Let $p = x^{m_1} (yx)^{(D - 4)/2}$.
We prove that
$W = W(x,y; m_1, 1, 1, \ldots,1, m_D)$ associates
by showing that $\pi^{D-2}(W) = \pi^{D-1}(W)$; that is,
$p y \cdot x y^{m_D} = p y x \cdot y^{m_D}$.
By Definition \ref{def:ARIF} applied twice, we have
$R(aba)R(b)R(a) = R(a)R(b)R(aba) $.
In particular,
if $a = y\iv$ and $b = yxy^{m_D + 1}$ then
$aba = xy^{m_D}$, so we get:
\[
R( xy^{m_D} )R( yxy^{m_D + 1})R(y\iv) =
R(y\iv)R( yxy^{m_D + 1})R( xy^{m_D} ) \ \ .
\]
We apply this equation to $p y^{-m_D} x\iv$:

Now, $ p y^{-m_D} x\iv\cdot xy^{m_D}$ 
is a product of a word with $D$ blocks, and by
Lemma \ref{lemma-dias-induct-cancel}, this is equal to
$p$ . Thus, applying Lemma \ref{lemma-dias-two} and power alternativity:
\begin{eqnarray*}
&&(p y^{-m_D} x\iv ) R( xy^{m_D} )R( yxy^{m_D + 1})R(y\iv) =
( p \cdot yxy^{m_D + 1} ) y\iv = \\
&&( p yx \cdot y^{m_D + 1} ) y\iv =
p yx \cdot y^{m_D} \ \ .
\end{eqnarray*}

Likewise, $ p y^{-m_D} x\iv \cdot y\iv$ is a product of a word with $D$ blocks,
of form $W(x,y; m_1, 1, \cdots ,1, -m_D , -1, -1)$.
This word associates by $H(D-1)$, since it is the same as
$W(x,y\iv ; m_1, -1, \cdots ,1, m_D , -1, 1)$. 
Thus, 
\[
p y^{-m_D} x\iv \cdot y\iv =
p \cdot y^{-m_D} x\iv y\iv =
p y \cdot y^{-m_D -1} x\iv y\iv \ \ ;
\]
the second ``$=$'' is obtained by applying Lemma \ref{lemma-powers},
with $a = y$ and $q = x\iv y\iv$.
We thus have
\begin{eqnarray*}
( p y^{-m_D} x\iv) R(y\iv)R( yxy^{m_D + 1})R( xy^{m_D} ) &=& \\
\left[(py \cdot y^{-m_D-1} x\iv y\iv )( yxy^{m_D + 1})\right]( xy^{m_D} ) &=&
py \cdot ( xy^{m_D} ) \ \ ,
\end{eqnarray*}
and hence $H^+(1)$.
\end{proof}

Finally, Lemmas \ref{lemma-dias-odd} and \ref{lemma-dias-even} complete
the proof of Theorem \ref{thm-dias}.

\section{Examples}
\label{sec-ex}

\begin{figure}[htb]
\begin{center}
{\includegraphics[scale = 0.5]{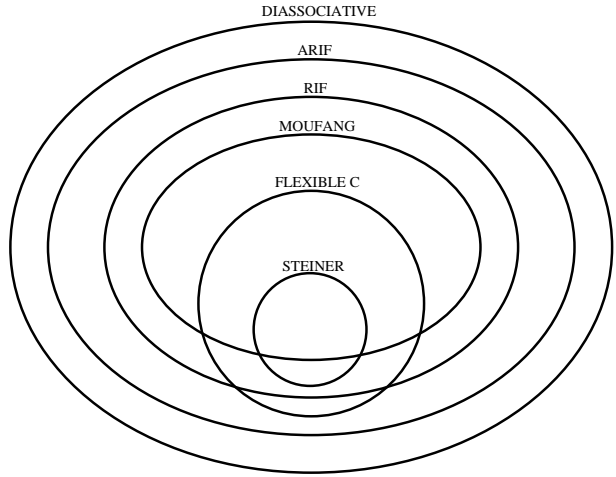}}
\end{center}
\caption{Some Varieties of Loops}
\label{fig-var}
\end{figure}

Figure \ref{fig-var} depicts the sub-varieties of diassociative
loops discussed in this paper.
All claimed inclusions have already been proved.
All regions shown are non-empty, as can easily be inferred
from results in the literature plus Example \ref{ex-flex-C-non-RIF}:
The loops which are both Moufang and Steiner are the boolean
groups and clearly are a proper sub-variety of the
extra loops, which are the loops which are both
Moufang and (flexible) C (see Fenyves \cite{FEN}),
and these in turn are properly contained in the Moufang loops.
If $A$ is, say, the 10-element Steiner loop, then it is not
a group and hence not Moufang.
The product of $A$ and any
non-boolean group will be a RIF flexible C-loop which is not Moufang
and not Steiner. 
The product of $A$ and any Moufang loop which is not an extra
loop will be a RIF loop which is not a C-loop.
Example \ref{ex-flex-C-non-RIF} is a flexible C-loop which is
not a RIF loop. Crossing this with a non-extra Moufang loop
yields a ARIF loop which is neither C nor RIF.
Finally, for every odd prime $p$, there is a diassociative
loop of order $p^3$ which is not a group; 
see, e.g., the proof of Theorem 5.2 in \cite{HK}.
Such loops cannot be Moufang by Chein \cite{CH}, and hence 
not ARIF by Corollary \ref{cor-odd}.

\begin{example}
\label{ex-flex-C-non-RIF}
There is a flexible C-loop which is not a RIF loop.
\end{example}
\begin{proof}
Consider the loop in Table \ref{table-flexC}.
The nucleus is $N = \{0,1,2\}$, and all squares are in $N$,
so that $L/N$ is the 8-element boolean group.
This is not RIF because
$(3 \cdot 12) \cdot (15 \cdot (3 \cdot 12))
\ne (3 \cdot ((12 \cdot 15) \cdot 3)) \cdot 12$,
so that (3) of Lemma \ref{lem:RIF-equiv} fails.
\end{proof}

\begin{table}[htb]
{ 
\footnotesize
\arraycolsep=1.2pt
\newcommand\bl{\bullet}
\[
\begin{array}{c|cccccccccccccccccccccccc|}
\bl& 0& 1& 2& 3& 4& 5& 6& 7& 8& 9&10&11&12&13&14&15&16&17&18&19&20&21&22&23 \\
\hline
0 & 0& 1& 2& 3& 4& 5& 6& 7& 8& 9&10&11&12&13&14&15&16&17&18&19&20&21&22&23 \\
1 & 1& 2& 0& 4& 5& 3& 7& 8& 6&10&11& 9&13&14&12&16&17&15&19&20&18&22&23&21 \\
2 & 2& 0& 1& 5& 3& 4& 8& 6& 7&11& 9&10&14&12&13&17&15&16&20&18&19&23&21&22 \\
3 & 3& 4& 5& 0& 1& 2& 9&10&11& 6& 7& 8&18&19&20&21&22&23&12&13&14&15&16&17 \\
4 & 4& 5& 3& 1& 2& 0&10&11& 9& 7& 8& 6&19&20&18&22&23&21&13&14&12&16&17&15 \\
5 & 5& 3& 4& 2& 0& 1&11& 9&10& 8& 6& 7&20&18&19&23&21&22&14&12&13&17&15&16 \\
6 & 6& 7& 8& 9&10&11& 0& 1& 2& 3& 4& 5&15&16&17&12&13&14&21&22&23&18&19&20 \\
7 & 7& 8& 6&10&11& 9& 1& 2& 0& 4& 5& 3&16&17&15&13&14&12&22&23&21&19&20&18 \\
8 & 8& 6& 7&11& 9&10& 2& 0& 1& 5& 3& 4&17&15&16&14&12&13&23&21&22&20&18&19 \\
9 & 9&10&11& 6& 7& 8& 3& 4& 5& 0& 1& 2&21&22&23&19&20&18&17&15&16&12&13&14 \\
10 & 10&11& 9& 7& 8& 6& 4& 5& 3& 1& 2& 0&22&23&21&20&18&19&15&16&17&13&14&12 \\
11 & 11& 9&10& 8& 6& 7& 5& 3& 4& 2& 0& 1&23&21&22&18&19&20&16&17&15&14&12&13 \\
12 & 12&14&13&18&20&19&15&17&16&21&23&22& 0& 2& 1& 6& 8& 7& 3& 5& 4& 9&11&10 \\
13 & 13&12&14&19&18&20&16&15&17&22&21&23& 1& 0& 2& 7& 6& 8& 4& 3& 5&10& 9&11 \\
14 & 14&13&12&20&19&18&17&16&15&23&22&21& 2& 1& 0& 8& 7& 6& 5& 4& 3&11&10& 9 \\
15 & 15&17&16&21&23&22&12&14&13&19&18&20& 6& 8& 7& 0& 2& 1&10& 9&11& 3& 5& 4 \\
16 & 16&15&17&22&21&23&13&12&14&20&19&18& 7& 6& 8& 1& 0& 2&11&10& 9& 4& 3& 5 \\
17 & 17&16&15&23&22&21&14&13&12&18&20&19& 8& 7& 6& 2& 1& 0& 9&11&10& 5& 4& 3 \\
18 & 18&20&19&12&14&13&21&23&22&17&16&15& 3& 5& 4&11&10& 9& 0& 2& 1& 6& 8& 7 \\
19 & 19&18&20&13&12&14&22&21&23&15&17&16& 4& 3& 5& 9&11&10& 1& 0& 2& 7& 6& 8 \\
20 & 20&19&18&14&13&12&23&22&21&16&15&17& 5& 4& 3&10& 9&11& 2& 1& 0& 8& 7& 6 \\
21 & 21&23&22&15&17&16&18&20&19&12&14&13& 9&11&10& 3& 5& 4& 6& 8& 7& 0& 2& 1 \\
22 & 22&21&23&16&15&17&19&18&20&13&12&14&10& 9&11& 4& 3& 5& 7& 6& 8& 1& 0& 2 \\
23 & 23&22&21&17&16&15&20&19&18&14&13&12&11&10& 9& 5& 4& 3& 8& 7& 6& 2& 1& 0 \\
\hline
\end{array}
\]
} 
\caption{A Flexible C non-RIF Loop}
\label{table-flexC}
\end{table}

\begin{example}
\label{ex-C-non-flex}
There is a C-loop which is not flexible.
\end{example}
\begin{proof}
Consider the loop in Table \ref{table-nonflexC}.
The nucleus is $N = \{0,1,2\}$, and all squares are in $N$,
so that $L/N$ is the 4-element boolean group.
This is not flexible because
$3 \cdot (6 \cdot 3) \ne (3 \cdot 6) \cdot 3$.
\end{proof}

\begin{table}[htb]
{ 
\footnotesize
\arraycolsep=1.2pt
\[
\begin{array}{c|cccccccccccc|}
\bullet & 0& 1& 2& 3& 4& 5& 6& 7& 8& 9&10&11 \\
\hline
0 & 0& 1& 2& 3& 4& 5& 6& 7& 8& 9&10&11 \\
1 & 1& 2& 0& 4& 5& 3& 7& 8& 6&10&11& 9 \\
2 & 2& 0& 1& 5& 3& 4& 8& 6& 7&11& 9&10 \\
3 & 3& 4& 5& 0& 1& 2&10&11& 9& 8& 6& 7 \\
4 & 4& 5& 3& 1& 2& 0&11& 9&10& 6& 7& 8 \\
5 & 5& 3& 4& 2& 0& 1& 9&10&11& 7& 8& 6 \\
6 & 6& 7& 8&11& 9&10& 0& 1& 2& 4& 5& 3 \\
7 & 7& 8& 6& 9&10&11& 1& 2& 0& 5& 3& 4 \\
8 & 8& 6& 7&10&11& 9& 2& 0& 1& 3& 4& 5 \\
9 & 9&10&11& 7& 8& 6& 5& 3& 4& 0& 1& 2 \\
10 & 10&11& 9& 8& 6& 7& 3& 4& 5& 1& 2& 0 \\
11 & 11& 9&10& 6& 7& 8& 4& 5& 3& 2& 0& 1 \\
\hline
\end{array}
\]
} 
\caption{A non-Flexible C-Loop}
\label{table-nonflexC}
\end{table}

The three examples in this section were found using 
the program SEM \cite{ZZ}, which simply outputs tables, such as
Tables \ref{table-flexC} and \ref{table-nonflexC}.
We do not see a really simple way of checking that 
Examples \ref{ex-flex-C-non-RIF} and \ref{ex-C-non-flex} are
both C-loops, with the first one also flexible.
However, the reader can easily write
the obvious computer program (entering each loop as an array)
to check these facts;
it is not necessary to verify that the code for SEM itself is correct.
Likewise, a program easily checks that the nucleus is $\{0,1,2\}$ for
both loops. On the other hand,
Example \ref{ex-steiner} below is a Steiner loop.
For this one, it was easy enough to match SEM's table
with the known triple systems, and then verify its properties
directly from facts about such systems.

The proof of diassociativity of ARIF loops in Section
\ref{sec-dias} is by induction on the number of blocks,
as is Moufang's proof for Moufang loops in \cite{M1,M2}, but
her proof is quite a bit shorter than ours.
She first shows that whenever $(vu)w = v(uw)$,
the same equation holds if the elements $u,v,w$ are permuted
or replaced by their inverses (\cite{M2}, pp.~420-421).
Using this fact, 
the step from 3 -- associativity to full diassociativity
is quite easy (the details are in \cite{M1}\S 1).
Actually, as Bruck pointed out, by using this fact
one can give a somewhat simpler ``maximal associative set'' argument
which avoids mentioning blocks at all (see \cite{BR}, \S VII.4).
However, as the following example shows,
this fact does not hold in all ARIF loops,
or even in all Steiner loops:

\begin{example}
\label{ex-steiner}
There is a Steiner loop of order 14 with elements
$u,v,w$ such that $(vu)w = v(uw)$ but $(uv)w \ne u(vw)$.
\end{example}
\begin{proof}
Let $L = \ZZZ_{13} \cup \{e\}$. Here, $e$ is the identity element
of the loop, so $x e = ex = x$ and $xx = e$ by definition.
Products $xy$ for
distinct elements $x,y$ of
$\ZZZ_{13} = \{0,1, \ldots 12\}$ are computed in the usual
way from a Steiner triple system $S$ on $\ZZZ_{13}$;
that is, $S$ is a set of 3-element subsets of $\ZZZ_{13}$, and
$xy = yx = z$, where $z$ is the (unique) element of
$\ZZZ_{13}$ such that $\{x,y,z\} \in S$.

For $S$, we take one of the standard examples of a triple
system (see, e.g., Example 19.12 of \cite{LW}):
$S$ contains blocks of the form $A_n = \{n, n+2, n+8\}$
and $B_n = \{n, n+3, n+4\}$, where $n \in \ZZZ_{13}$.

So, for example
$1 \cdot 0 = 10$ (since $B_{10} = \{10,0,1\}$), 
$10 \cdot 5 = 12$ (using $A_{10}$),
$0 \cdot 5 = 7$ (using $A_5$,), and $1 \cdot 7 = 12$ (using $A_{12}$).
Thus, $(1 \cdot 0) \cdot 5 = 1 \cdot (0 \cdot 5) = 12$.

However, $(0 \cdot 1) \cdot 5 = 12 \ne 3 = 0 \cdot (1 \cdot 5)$,
since $1 \cdot 5 = 4$ (using $B_1$) and $0 \cdot 4 = 3$ (using $B_0$).
\end{proof}

\bigskip
\noindent{\textit{Acknowledgement:}} Thanks to the referee for
useful suggestions.

\end{document}